\documentclass[12pt,oneside]{amsart}

\usepackage{amssymb}
\usepackage{rotating}

\usepackage{amsfonts}
\usepackage{dsfont}
\usepackage{amsxtra}

\usepackage{index}

\usepackage{amsmath}

\usepackage{index}

\theoremstyle{plain}
\newtheorem{theorem}{Theorem}[section]
\newtheorem{lemma}[theorem]{Lemma}
\newtheorem{proposition}[theorem]{Proposition}
\newtheorem{corollary}[theorem]{Corollary}

\theoremstyle{definition}
\newtheorem{definition}[theorem]{Definition}

\theoremstyle{remark}
\newtheorem{remark}[theorem]{Remark}

 \numberwithin{equation}{subsection}

 \newenvironment{dem}{\noindent{\it Proof}.}{\qed}

\newcommand{\ff}{\mathbb{F}}

 \newcommand{\ext}{\mathrm{Ext}}

      \makeatletter
      \def\@setcopyright{}
      \def\serieslogo@{}
      \makeatother

\begin{document}

\author{Julio Jos\'e Moyano-Fern\'andez}
\address{Institut f\"ur Mathematik, Universit\"at Osnabr\"uck. Albrechtstra\ss e 28a, D-49076 Osnabr\"uck, Germany}
\email{jmoyanof@uni-osnabrueck.de}

\title[On Weierstra\ss~ semigroups and Poincar\'e series]{On Weierstra\ss~ semigroups at one and two points and their corresponding \\ Poincar\'e series}

\begin{abstract}
The aim of this paper is to introduce and investigate the Poincar\'e
series associated with the Weierstra\ss~semigroup of one and two
rational points at a (not necessarily irreducible) non--singular
projective algebraic curve defined over a finite field, as well as
to describe their functional equations in the case of an affine
complete intersection.
\end{abstract}



 \thanks{Supported partially by the grant of the Spanish Government
``Ministerio de Educaci\'on MTM2007--64704", in cooperation with the European Union in the framework of the founds ``FEDER'', by the grant of the Deutsches Akademischen
Austauschdienst (DAAD)--La Caixa, and by the Deutsche
Forschungsgemeinschaft (DFG)}

\maketitle

\section{Introduction}

The origin of the Weierstra\ss~ semigroup lies in the
``L\"uckensatz'' pointed out by Karl Weierstra\ss~ in some of his
lectures in the 1860's: For every point $P$ on a compact Riemann
surface $X$ of genus $g$, there are exactly $g$ integers $l_i(P)$
with
$$
1=l_1 (P) < \ldots < l_g(P)=2g-1
$$
so that there is no meromorphic function on $X$ having a pole at
$P$ of multiplicity $l_i (P)$ as its only singularity (see
\cite[III, pp.~297--307]{weier}). The set $G(P):=\{l_i (P)\}$ is
called the gap set of $P$. After some remarkable results of Brill
and Noether (\cite{BN}), Hurwitz realised that, if $\alpha,
\beta$ belong to $\mathbb{N} \setminus G(P)$ and $z_{\alpha},z_{\beta}$
are two functions having a pole of order $\alpha$ resp. $\beta$ at
$P$ as their only singularities, then the function $z_{\alpha}
z_{\beta}$ has a pole of order $\alpha + \beta$ at $P$ as its only
singularity (cf. \cite{cen}). It implies that $\mathbb{N} \setminus G(P)$ is a
semigroup (cf. \cite{hur}, p.~409 after equation (6)): The
Weierstra\ss~semigroup associated with the point $P$.
\medskip

Let $X$ be a non--singular projective algebraic curve defined over
a finite field $\ff$. Let $P_1, \ldots ,P_r$ be a set of rational
points on $X$ and consider the family of (finite dimensional)
$\ff$--vector subspaces 
\[
\Gamma (X,\underline{mP}):=\Gamma(X,m_1
P_1 + \ldots + m_r P_r)
\]
with $m_i \in \mathbb{Z}$ for all $i=1,
\ldots , r$ (cf. Subsect. \ref{sec:rr} below). This family gives
rise to a $\mathbb{Z}^r$--multi--index filtration on the
$\ff$--algebra $A:= \bigcup_{\underline{m} \in \mathbb{Z}^r}
\Gamma (X, \underline{mP})$ of the affine curve $X \setminus
\{P_1, \ldots , P_r\}$. This multi--index filtration is related to
Weierstra\ss~semigroups (with respect to several points in
general, see \cite{delgado}) and, in the case of finite fields, to
the methodology for trying to improve the Goppa estimation of the
minimal distance of algebraic--geometrical codes (see
\cite{cato}). Furthermore,  a connection of that
filtration with global geometrical--topological aspects in a
particular case  is shown in \cite{cadegu}. \medskip

Poincar\'e series are typically series
for which the coefficients represent discrete information about
the class of the study objects---the Weierstra\ss~semigroup in
this case. Thus, a natural question, which this work is devoted to, is to
define in a proper way the notion of Poincar\'e series and to
describe it, in the case of a Weierstra\ss~ semigroup associated
with one or two points. This will be made by associating such a series to 
the dimensions of the vector spaces $\Gamma (X,\underline{mP})$.
\medskip

The article goes as follows. Sect.~\ref{sec:1} is devoted to the definition and
main properties of the Weierstra\ss~semigroup for not necessarily
irreducible curves using the notion of Manis valuation. In Sect.~\ref{sec:2} 
we introduce and describe the Poincar\'e series of the
Weierstra\ss~semigroup associated with one and two points.
A formula for the conductor of the semigroup for the case $X$ to 
be non--singular in the affine part, as well as descriptions for the 
genus of $X$ and the gaps of the semigroup if the curve is affine
complete intersection are provided. Finally we
consider symmetry properties in Sect.~\ref{sec:3} (revising the
work of Delgado \cite{delgado}), from which we deduce functional
equations for the Poincar\'e series.

\section{Weierstra\ss~semigroup for reducible curves}
\label{sec:1}

\subsection{Reducible curves and Manis valuations}\label{sec:manis}

We can extend the concept of valuation to start from a ring
satisfying certain properties, instead of a field. More precisely,
let $\mathcal{K}$ be a ring with large Jacobson radical (i.e., a
ring in which every prime ideal containing the Jacobson radical is
maximal; for instance, a semilocal ring) and being its own total
ring of fractions. A surjective map $w: \mathcal{K} \to
\mathbb{Z}_{\infty}:=\mathbb{Z} \cup \{ \infty \}$ such that $w(1)=0$,
$w(0)=\infty$, $w(ab)=w(a)+w(b)$, for all $a,b \in \mathcal{K}$,
and $w(a+b) \ge \min \{w(a),w(b) \}$, for all $a,b \in \mathcal{K}$,
is called a \emph{Manis valuation} of the ring $\mathcal{K}$ (cf.
\cite[Chapter I, Def. 2.8, p.~9]{kiyek}).
\medskip

Let $X$ be a reducible projective curve of arithmetical genus $g_X$
defined over a finite field $\ff$, which is locally a complete
intersection (or simply a \emph{curve} from now on). Let $X=X_1
\cup \ldots \cup X_r$ be its decomposition into irreducible
components.  A function field $K(X_i)$ can be associated to each $X_i$,
and $K(X):=K(X_1) \times \ldots \times K(X_r)$ is the ring of
rational functions over $X$. A point of $X$ is called
\emph{regular} if it lies on just one irreducible component of $X$
and if it is a simple point of this component. Let us take $P$ a
regular point of $X$. Then there exists some $\kappa \in \{1,
\ldots , r \}$ such that $P \in X_{\kappa}$ and the local ring of
$X$ at $P$ coincides with the local ring of $X_{\kappa}$ at $P$,
and therefore, one can associate with $P$ a discrete valuation
$v_P:K(X_{\kappa}) \to \mathbb{Z}_{\infty}$. Let us consider
$z=(z_1, \ldots, z_r) \in K(X)$. For every $1 \le j \le r$, we
define
\begin{displaymath}
w_P^{\prime}(z_j)=\left\{%
\begin{array}{ll}
    v_P(z_j), & \hbox{~if $j=\kappa$;} \\
    \infty, & \hbox{~otherwise.} \\
\end{array}%
\right.
\end{displaymath}

It is easy to check that the map $w_P: K(X) \to
\mathbb{Z}_{\infty}$, $z \mapsto \min_{1 \le j \le r} w^{\prime}_P
(z_j)$ is a discrete Manis valuation of the ring $K(X)$.
\medskip

Let $\underline{P}:= \{ P_1, \ldots , P_r\}$ be a set of non--singular
points on $X$, and set $w_i=w_{P_i}$ a discrete Manis valuation of
$K(X)$ associated with each $P_i$, for all $1 \le i \le r$. If we
consider the affine coordinate ring
$\Gamma(X^{\prime},\mathcal{O}_X)$ of the affine curve
$X^{\prime}:=X \setminus \underline{P}$, then we may define the
Weierstra\ss~semigroup of $X$ at $\underline{P}$ as the following
sub--semigroup of $(\mathbb{Z}_{\infty})^{r}$:
\[
\Gamma_{\underline{P}}:=\left \{ -(w_1,\ldots,w_r)(f) \mid f \in
\Gamma(X^{\prime},\mathcal{O}_X) \right \}.
\]


\emph{Notations}. Let $\underline{n}=(n_1, \ldots, n_r) \in
\mathbb{Z}^r$, $J=\{i_1,\ldots,i_d \} \subseteq I:=\{1, \ldots,r \}$.
Consider the subsets
\[
\overline{\nabla}_J (\underline{n})= \{(m_1, \ldots,m_r) \in
\mathbb{Z}^r \mid m_i=n_i~ \forall i \in J,~ m_j < n_j~ \forall j
\notin J \};
\]
\[
\overline{\nabla}(\underline{n}) = \bigcup_{i=1}^{r}
\overline{\nabla}_i (\underline{n});
\]
and $\nabla_J (\underline{n}):=\overline{\nabla}_J (\underline{n})
\cap \Gamma_{\underline{P}}$, $\nabla(\underline{n}):=
\overline{\nabla}(\underline{n}) \cap \Gamma_{\underline{P}}$, as
well as $\nabla_i^r (\underline{n}):= \bigcup_{i \in J} \nabla_J
(\underline{n})$.
The vectors of the standard basis of $\mathbb{Z}^r$ will be denoted by $e_1, \ldots, e_r$.
\medskip

\begin{definition}
An element $\underline{n} = (n_1, \ldots , n_r)$ of the semigroup $\Gamma_{\underline{P}}$
is said to be maximal for $\Gamma_{\underline{P}}$ if $\nabla(\underline{n})=\varnothing$. If, \emph{moreover}, $\nabla_J(\underline{n})=\varnothing$ for every $J \subsetneq I$, then the maximal element is called absolute maximal. A maximal element  for $\Gamma_{\underline{P}}$ is called relative maximal if $\nabla_J (\underline{n}) \ne \varnothing$ for every $J \subsetneq I$ with $\sharp J \ge 2$ 
(and $\nabla_J (\underline{n}) = \varnothing$ if $\sharp J = 1$).
\end{definition}

\begin{remark}
For the case $r=2$, the concepts of relative and absolute maximal
coincide. Thus we will simply say ``maximal'' points of
$\Gamma_{P_1,P_2}$. The set of maximal points of
$\Gamma_{P_1,P_2}$ will be denoted by $\mathcal{M}_{P_1,P_2}$.
\end{remark}

\subsection{The Riemann--Roch theorem for reducible curves}
\label{sec:rr}

Let $X$ be a curve. Let $\mathcal{O}_X$ be the sheaf of local
rings on $X$. We can take, for each invertible sheaf $\mathcal{L}$
on $X$, a divisor $D$ with support contained in the set of regular
points of the curve such that $\mathcal{L}=\mathcal{L}(D)$, where
$\mathcal{L}(D)$ is the divisorial sheaf related to $D$. As
$\mathcal{L}(D)$ is coherent, the cohomology groups
$H^q(X,\mathcal{L}(D))$ are $\ff$--vector spaces of finite
dimension (if $q \ge 2$ its dimension is $0$). Furthermore, since
$\mathcal{L}(D) \cap \mathcal{O}_X$ is a coherent sheaf too, we
can define the \emph{degree} of the divisor $D$ as
\[
\deg (D) = \dim_{\ff} \Gamma \left ( X,
\mathcal{L}(D)/(\mathcal{L}(D) \cap \mathcal{O}_X) \right ) -
\dim_{\ff} \Gamma \left (X, \mathcal{O}_X / (\mathcal{L}(D) \cap
\mathcal{O}_X) \right).
\]
Such a definition extends in a natural way the classical one. Since $X$ is locally a complete intersection, the dualizing sheaf of
$X$ is an inver\-tible sheaf on $X$ and we can take the
corresponding canonical divisor $K$ such that $\mathrm{supp}(K)$
is contained in the set of regular points on $X$. By setting
\begin{displaymath}
\begin{array}{lllll}
\ell(D) & = & \dim_{\ff} H^0 (X, \mathcal{L}(D))& = & \dim_{\ff} \Gamma (X, \mathcal{L}(D))\\
i(D) & = & \dim_{\ff} H^1 (X, \mathcal{L}(D)) & = & \ell(D - K)\\
\gamma & = & \dim_{\ff} H^0 (X,\mathcal{O}_X) & = & \mathrm{number
~ of ~ connected~ parts~ of~} X,
\end{array}
\end{displaymath}
the Riemann--Roch theorem for reducible curves can be proven
(see~\cite[p.~103]{oort}):
\begin{theorem}[Riemann--Roch] \label{thm:RR}
For every divisor $D$ on $X$, we have
\[
\ell(D)-i(D)=\deg (D) - g_X + \gamma.
\]
\end{theorem}

If $D= \{(U_i,f_i)\}_{i \in I}$ is a Cartier divisor in $X$, and
$P$ is a regular rational point on $X$, we denote by $w_P(D)$ the
value $w_P(f_j)$ if $P \in U_j$ (in the sense of Manis). Usual
notations for effectivity of $D$ (namely $D > 0$) and linear
equivalence between two divisors ($D \sim D^{\prime}$) hold.
\medskip

Let us consider a set $\{P_1, \ldots, P_r\}$ of non--singular rational
points on $X$, a vector $\underline{n}=(n_1, \ldots, n_r) \in \mathbb{Z}^r$
and the divisor $\underline{nP}=\sum_{i=1}^{r}n_i P_i$. Denote the
$\ff$--vector space of global sections of the invertible sheaf
$\mathcal{L}(\underline{nP})$ by $\Gamma (X,\underline{nP})$, its
dimension by $\ell(\underline{nP})$ and $\ell (K -
\underline{nP}):=\dim_{\ff} \Gamma (X,
\mathcal{L}(K-\underline{nP}))$. We have
\begin{eqnarray}
\Gamma (X,\underline{nP})&=& \{ f \in K(X)^{\ast} \mid (f) +
\underline{nP} > 0 \} \nonumber \\
&=& \{f \in K(X)^{\ast} \mid w_i(P) \ge -n_i, ~ i=1, \ldots,r
\},\nonumber
\end{eqnarray}
where $K(X)^{\ast}:=K(X) \setminus \{0\}$. Notice that $\Gamma (X,
\underline{nP}) \supseteq \Gamma (X,\underline{mP})$ for every
$\underline{n}, \underline{m} \in \mathbb{Z}^r$ with $n_i \ge
m_i$, for all $i \in \{1, \ldots ,r\}$; hence the set $\{\Gamma
(X, \underline{nP})\}_{\underline{n} \in \mathbb{Z}^r}$ defines a
filtration given by multi--indices (or multi--index filtration) of
the curve.

\subsection{Main properties of the semigroup}

We collect now some basic properties and definitions around the
Weierstra\ss~semigroup.

\begin{definition}
An element $\underline{m} \in \mathbb{Z}^r$ is said to be a gap of
$\Gamma_{\underline{P}}$ if $\underline{m} \notin \Gamma_{P}$.
Otherwise $\underline{m}$ is said to be a non-gap of
$\Gamma_{\underline{P}}$.
\end{definition}

Let us mention two important properties of these concepts 
(their proofs, also valid for Manis valuations, can be found in \cite{delgado}; or \cite[Sect.~2]{cato}).

\begin{lemma}
Let $\underline{m} \in \mathbb{Z}^r$, then $\underline{m}$ is a
non--gap of $\Gamma_{\underline{P}}$ if and only if for all $i \in
\{1, \ldots, r\}$ one has $\ell (\underline{mP})= \ell ((\underline{m}-e_i)P)+1$.
\end{lemma}

\begin{lemma}
If $X$ is a curve of genus $g$ and $\underline{m} = (m_1, \ldots ,
m_r ) \in \mathbb{Z}^r$ is a gap, then $0 \le m_1 + \ldots + m_r
<2g_X$.
\end{lemma}

\begin{definition}
Let $\tau \in \mathbb{Z}^r_{> 0}$. The Weierstra\ss~semigroup $\Gamma_{\underline{P}}$ is said to be symmetric with respect to $\tau$ if it satisfies the following property (cf. \cite{delgado}):
``If $\underline{n} \in \mathbb{Z}^r$, then $\underline{n} \in \Gamma_{\underline{P}}$ if and only if $\nabla(\tau-\underline{n})=\varnothing$".
\end{definition}

Notice that for $r=1$ this means: The semigroup $\Gamma_P$ is symmetric if there exists an integer $n$ such that the mapping $\mathbb{Z} \rightarrow \mathbb{Z}$, given by $z \mapsto n-z$, maps elements of $\Gamma_P$ onto gaps and gaps onto elements of $\Gamma_P$ (cf. definition in \cite[p.~32]{hk}). Furthermore, the case $r=1$ is remarkable in the following sense:

\begin{lemma}
Let $r=1$. Then $\Gamma_P$ is a numerical semigroup.
\end{lemma}

This implies that $\Gamma_P$ possesses a conductor, denoted by $c(\Gamma_P)$, when $r=1$ (see e.~g.~\cite{buch}). The gap $g(\Gamma_P):=c(\Gamma_P)-1$ is called the Frobenius number associated to the semigroup $\Gamma_P$.

\begin{lemma} \label{lem:period}
Let $X$ be a curve, $P_1, \ldots, P_r$ rational points on $X$ and
$\Gamma_{\underline{P}}$ the Weierstra\ss~semigroup of $P_1,
\ldots, P_r$. Let $\underline{m}:=(m_1, \ldots, m_r) \in
\mathbb{Z}^r$.
\begin{itemize}
  \item[(i)] There exists a number $\vartheta \in \mathbb{Z} \setminus \{0\}$ so that
  $\vartheta(P_i-P_j)$ is a principal divisor, for any $i, j \in \{1, \ldots r\}$ with $i \ne j$.
  \item[(ii)] For any $\lambda \in \mathbb{Z}$, then $\underline{m} \in
  \Gamma_{\underline{P}}$ if and only if $\underline{m} + \lambda (\vartheta,-\vartheta,0, \ldots , 0) \in
  \Gamma_{\underline{P}}$.
\end{itemize}
\end{lemma}

\dem~(i) Let $\mathrm{Pic}(X)$ be the set of isomorphism classes
of invertible sheaves on $X$, and $\mathrm{Pic}^0(X)$ the subgroup
of $\mathrm{Pic}(X)$ of those classes of degree $0$. The Jacobian
variety of $X$ is isomorphic to $\mathrm{Pic}^0(X)$ (cf.
\cite[Chapter 7, Theorem 4.39]{liu}). Since the
$\mathbb{F}$--points of $J$ are a finite group, then the divisor
$P_i-P_j \in \mathrm{Pic}^0(X)$ is an element of finite order, for
any $i,j \in \{1, \ldots , r\}$, $i \ne j$. (ii) is a consequence
of (i). \qed
\medskip

We are in this article much more interested in the case $r=2$. Let us assume 
for the rest of this section $\Gamma_{\underline{P}}$ to be the Weierstra\ss~semigroup 
of two rational points $P_1,P_2$ on $X$.

\begin{definition}
Let us assume $\sharp(\ff) > 2$. The smallest natural number
$\vartheta \ne 0$ such that $\vartheta(P_1-P_2)$ is a principal
divisor is called the \emph{period} of the semigroup
$\Gamma_{P_1,P_2}$.
\end{definition}

\begin{remark}
Consider the discrete plane $\mathbb{Z}^2$. The period of the
semigroup defines an important area on it, namely the parallelogram in
$\Gamma_{P_1,P_2} \cap \mathbb{Z}^2$ determined by the vertices
$(0,0),(0,2g),(\vartheta,2g-\vartheta)$ and
$(\vartheta,-\vartheta)$ not including the line joining $(0,0)$
and $(0,2g)$, which will be called the \emph{fundamental corner} of
$\Gamma_{P_1,P_2}$, and will be denoted by
$\mathcal{C}_{P_1,P_2}=\mathcal{C}$. By Lemma \ref{lem:period}
(ii), any translation of $\mathcal{C}_{P_1,P_2}$ of vector
$\lambda(\vartheta,-\vartheta)$ for every $\lambda \in \mathbb{Z}
\setminus \{0\}$, reproduces the same distribution of points in
$\Gamma_{P_1,P_2}$. Therefore, the points in $\mathcal{C}$ will
play an essential role in the semigroup's knowledge, as we will
see in the next sections, especially the set of maximal points of
$\mathcal{C}$, i.e. the set
$\mathcal{M}_{\mathcal{C}}:=\mathcal{M}_{P_1,P_2} \cap \mathcal{C}$.
\end{remark}

\section{Poincar\'e series of the Weierstra\ss~semigroups associated with one and two
points}\label{sec:2}

This section is devoted to introduce the notion of Poincar\'e series
associated with the multi--index filtration given by $\Gamma (X,
\underline{mP})$ and to show its behaviour in the case of one and
two points. For every $\underline{m} \in \mathbb{Z}^r$, set
\[
d(\underline{m}):=\dim_{\ff} \frac{\Gamma
(X,\underline{mP})}{\Gamma(X,(\underline{m}-\underline{1})\underline{P})}.
\]

\begin{definition} \label{defn:PS}
The Poincar\'e series of the Weierstra\ss~semigroup
$\Gamma_{\underline{P}}$ is defined to be
\[
P_{\Gamma_{\underline{P}}}(\underline{t}):= \sum_{\underline{m}
\in \mathbb{Z}^r} d(\underline{m}) \underline{t}^{\underline{m}},
\]
where $\underline{t}^{\underline{m}}:=t_1^{m_1} \cdot \ldots \cdot
t_r^{m_r} $.
\end{definition}

The series $P_{\Gamma_{\underline{P}}}(\underline{t})$ contains
negative powers of variables $t_i$, however it does not contain
monomials $\underline{t}^{\underline{m}}$ with purely negative
$\underline{m}$ (i.e., all components of which are negative). It
will be convenient to consider
$P_{\Gamma_{\underline{P}}}(\underline{t})$ as an element of the
set of formal Laurent series in the variables $t_1, \ldots, t_r$
with integer coefficients without purely negative exponents; i.e.,
of expressions of the form $\sum_{\underline{m} \in \mathbb{Z}^r
\setminus \mathbb{Z}^r_{\le -1}} d(\underline{m})
\underline{t}^{\underline{m}}$.

\subsection{Poincar\'e series associated with one point}

Consider $\Gamma_{P_1}=\Gamma_{P}$ the Weierstra\ss~semigroup
associated with one rational point on $X$. In this case we have that $d(m)=1$ if and
only if $m \in \Gamma_P$, hence $P_{\Gamma_P}(t)=\sum_{n \in
\Gamma_P} t^n$.
\medskip

Assume $X$ to be smooth on the affine part. Then $\Gamma_P$ is
strictly generated by some elements $\{r_0, \ldots , r_h\}$ (cf.
\cite{ab}). By setting 
\[
\theta_i:=\mathrm{gcd}(r_0, \ldots ,
r_{i-1})
\]
for $1 \le i \le h+1$ and  $d_i:= \theta_{i} /
\theta_{i+1}$ for $1 \le i \le h$, the semigroup $\Gamma_P$
satisfies the following property: Every element $n \in \Gamma_P$
can be written uniquely as
\[
n=\lambda_0 r_0 + \lambda_1 r_1 + \ldots + \lambda_h r_h, \eqno
(\ast)
\]
with $\lambda_0 \ge 0$ and $0 \le \lambda_i < d_i$, for $1 \le i
\le h$. This leads to the following description of the Poincar\'e
series $P_{\Gamma_P}(t)$ (cf. \cite[p.~1840]{cadegu}):

\begin{proposition} \label{prop:uno}
Let $X$ be non--singular in the affine part. Let $\Gamma_P$ be the
Weierstra\ss~semigroup at a point $P$ on $X$. We have
\[
P_{\Gamma_P}(t)= \frac{1}{1-t^{r_0}} \cdot \prod_{i=1}^{h}
\frac{1-t^{n_i r_i}}{1-t^{r_i}}=\frac{1}{1-t^{r_0}}
\prod_{j=1}^{h} \sum_{i=0}^{n_j-1} t^{i r_j}.
\]
In particular, the Poincar\'e series $P_{\Gamma_P}(t)$ is a rational
function.
\end{proposition}

\dem~ The formula can be proven by considering Equation $(\ast)$ as follows:
\begin{align*}
\sum_{n \in \Gamma_P} t^n & =
 \sum_{\substack{\lambda_0 \ge 0 \\ 0 \le \lambda_i < n_i \\ 1 \le i \le h}}
 t^{\lambda_0 r_0 + \ldots + \lambda_h r_h}    \\
     & =   \Big ( \sum_{\lambda_0 \ge 0} t^{\lambda_0 r_0} \Big )
      \cdot \Big ( \sum_{0 \le \lambda_1 < n_1} t^{\lambda_1 r_1} \Big )
       \cdot \ldots \cdot \Big ( \sum_{0 \le \lambda_g < n_h} t^{\lambda_h r_h}
        \Big )   \\
     & =  \frac{1}{1-t^{r_0}} \cdot \frac{1-t^{n_1
     r_1}}{1-t^{r_1}}\cdot \ldots \cdot \frac{1-t^{n_h
     r_h}}{1-t^{r_h}}.
\end{align*}
\qed

\begin{proposition} \label{prop:dos}
We have
\[
P_{\Gamma_P}(t)=\frac{L_{\Gamma_P} (t)}{1-t},
\]
with $L_{\Gamma_P}(t)=t^c +(1-t)\sum_{\substack{n \in \Gamma_P \\
n < c }} t^n \in \mathbb{Z}[t]$.
\end{proposition}

\begin{proposition}\label{prop:conductor}
Let $X$ be a non--singular curve in the affine part, let $P$ be a closed point on $X$. The conductor of $\Gamma_P$ is given by
\[
c(\Gamma_P)=1-r_0+\sum_{j=1}^{h} (n_j-1)r_j;
\]
therefore, the Frobenius number of $\Gamma_P$ is given by
\[
g(\Gamma_P)=\sum_{j=1}^{h} (n_j-1)r_j -r_0.
\]

\end{proposition}

\dem~By Proposition \ref{prop:uno} one has $(1-t^{r_0})P_{\Gamma_P}(t) = \widetilde{L}_{\Gamma_P}(t)$,
where
\[
\widetilde{L}_{\Gamma_P}(t)=\prod_{j=1}^{h}(1+t^{r_j}+t^{2r_j}+ \ldots + t^{(n_j-1)r_j}). 
\]
This together with Proposition \ref{prop:dos} leads to the equality
$(1+t+ \ldots +t^{r_0-1}) L_{\Gamma_P}(t)=\widetilde{L}_{\Gamma_P}(t)$. By comparing the leading terms on both sides of the equality we obtain the result.
\qed

\begin{corollary} 
Let $X$ be non--singular in the affine part and affine complete
intersection. Then the genus of $X$ is given by the equality
\[
g_X=\sum_{j=1}^{h} \frac{(n_j-1)r_j}{2} +\frac{1-r_0}{2}.
\]
\end{corollary}

\dem~ The equality is a consequence of the fact that $X$ is an
affine complete intersection if and only if $c=2g_X$ (cf.
\cite{sathaye}).\qed
\medskip

Another consequence of Proposition \ref{prop:conductor} is a characterisation of the gaps of $\Gamma_P$ for an
affine complete intersection curve.

\begin{corollary}
Let $X$ be non--singular in the affine part and affine complete
intersection. Let $m \in \mathbb{Z}$. Then $m$ is a gap of
$\Gamma_P$ if and only if there exist $\lambda_0, \ldots ,
\lambda_h > 0$ such that
\[
m=\sum_{i=1}^{h} n_i r_i - \sum_{j=0}^{h}\lambda_j r_j.
\]
\end{corollary}

\dem~ Assume $m \notin \Gamma_P$. Since the semigroup is symmetric under these hypothesis (see \cite{sathaye}), and by the formula for the computation 
of the Frobenius number of $\Gamma$ given in Proposition \ref{prop:conductor}, we have that
\[
\sum_{j=1}^{h} (n_j-1)r_j - r_0 -m \in \Gamma_P.
\]
Therefore, there exist $\alpha_0, \ldots , \alpha_h \ge 0$ such that
\[
\sum_{j=1}^{h} (n_j-1)r_j - r_0 -m = \sum_{i=0}^{h} \alpha_i r_i,
\]
and we deduce that
\[
m=\sum_{i=1}^{h} n_i r_i - \sum_{j=0}^{h}(\alpha_j +1) r_j.
\]
Conversely, suppose that $m=\sum_{i=1}^{h} n_i r_i - \sum_{j=0}^{h}\lambda_j r_j$ for $\lambda_j > 0$ 
for all $0 \le j \le h$.  Assuming that $m \in \Gamma_P$, then
\[
\sum_{i=1}^{h} n_i r_i - \sum_{j=0}^{h}\lambda_j r_j + \sum_{i=0}^{h} (\lambda_i -1) r_i \in \Gamma_P,
\]
but on the other hand this is in fact equal to the Frobenius number (by Proposition~\ref{prop:conductor}), i.e., 
a gap of $\Gamma_P$, which is a contradiction.
\qed

\subsection{Poincar\'e series associated with two points}

Let $\Gamma_{P_1,P_2}$ be the Weierstra\ss~semigroup associated
with two rational points $P_1, P_2$ on $X$. Let us define the
series
\[
L_{\Gamma_{P_1,P_2}}(\underline{t}):=(1-t_1)(1-t_2)P_{\Gamma_{P_1,P_2}}(\underline{t}).
\]
This expression can be rewritten as
\[
L_{\Gamma_{P_1,P_2}}(\underline{t})=\sum_{\underline{m} \in
\mathbb{Z}^2} c(\underline{m})\underline{t}^{\underline{m}},
\]
where
$c(\underline{m}):=d(\underline{m})-d(\underline{m}-e_1)-d(\underline{m}-e_2)+d(\underline{m}-\underline{1})$.
We will describe the series $L_{\Gamma_{P_1,P_2}}(\underline{t})$
in terms of the set $\mathcal{M}_{P_1,P_2}$ of maximal points of
$\Gamma_{P_1,P_2}$ in Proposition \ref{prop:clave}. First, a summary of 
the most elementary properties of the dimensions
$d(\underline{m})$, for $\underline{m}=(m_1,m_2) \in \mathbb{Z}^2$ is
presented. This can be found in Carvalho--Torres \cite[Sect.~2]{cato}

\begin{lemma} \label{lem:propertiessemigroup}
Let $\Gamma_{P_1,P_2}$ be the Weierstra\ss~semigroup associated
with two rational points $P_1, P_2$ on $X$. The following statements hold:
\begin{enumerate}
  \item $\underline{m} \in \Gamma_{P_1,P_2}$ if and only if $d(\underline{m}-e_i)
  =1$ for $i \in \{1,2\}$.
  \item Let $\underline{m} \in \Gamma_{P_1,P_2}$. Then $d(\underline{m})=1$ if
  and only if $\nabla (\underline{m}) = \varnothing$.
  \item $d(\underline{m})=2$ if and only if $\underline{m} \in \Gamma_{P_1,P_2}$ and $\nabla (\underline{m})
  \ne \varnothing$.
  \item If $m_1 \in \Gamma_{P_1}$ and $m_2 \in \Gamma_{P_2}$, then
  $d(\underline{m})=2$.
\end{enumerate}
\end{lemma}

\begin{proposition} \label{prop:clave}
Let $\underline{m} \in \mathbb{Z}^2$. Then we have
\begin{enumerate}
    \item $c(\underline{m})=-1$ if and only
if $\underline{m}-\underline{1}$ is maximal.
    \item $c(\underline{m})=1$ if and only if $\underline{m}$ is
    maximal.
\end{enumerate}
\end{proposition}

\dem~ Taking into account the main properties of the
Weierstra\ss~semigroup of Lemma \ref{lem:propertiessemigroup}, and
since
$c(\underline{m})=d(\underline{m})-d(\underline{m}-e_1)-d(\underline{m}-e_2)+d(\underline{m}-\underline{1})$,
we consider the different values the coefficients
$d(\underline{m})$ can take. So, if $d(\underline{m})=2$, then the
following statements hold:
\begin{itemize}
  \item[--] if $d(\underline{m}-e_1)=d(\underline{m}-e_2)=1$, then the only possible case is
  $d(\underline{m}-\underline{1})=0$ and
  $c(\underline{m})=2-1-1+0=0$;
  \item[--] if $d(\underline{m}-e_1)=1$ and $d(\underline{m}-e_2)=2$, then $d(\underline{m}-\underline{1})=1$ and from this
  $c(\underline{m})=2-1-2+1=0$;
  \item[--] if $d(\underline{m}-e_1)=2$ and $d(\underline{m}-e_2)=1$, as in the previous case
  one has $d(\underline{m}-\underline{1})=1$ and
  $c(\underline{m})=2-2-1+1=0$;
  \item[--] if $d(\underline{m}-e_1)= d(\underline{m}-e_2)=2$, one has $d(\underline{m}-\underline{1}) \ge 1$. If
  $d(\underline{m}-\underline{1})=2$, then $c(\underline{m})=2-2-2+2=0$; on the other hand, if
  $d(\underline{m}-\underline{1})=1$, then $c(\underline{m})=-1$ and $\nabla (\underline{m}-\underline{1}) = \varnothing$ by Lemma \ref{lem:propertiessemigroup} (2).
  Thus the point $\underline{m}-\underline{1}$ is maximal.
\end{itemize}

In the case $d(\underline{m})=1$, we have to distinguish whether
the point $\underline{m}$ is maximal. If this is the case, 
\begin{itemize}
  \item[--] if $d(\underline{m}-e_1) = d(\underline{m}-e_2)=1$, one
  has either (i) $d(\underline{m}-\underline{1})=1$ implying that $c(\underline{m})=1-1-1+1=0$ and
  $\underline{m}-\underline{1}$ is maximal; or (ii) $d(\underline{m}-\underline{1})=2$ with
  $c(\underline{m})=1-1-1+2=1$, hence $\underline{m}-\underline{1}$ is not
  maximal by Lemma \ref{lem:propertiessemigroup}(3);
  \item[--] if $d(\underline{m}-e_1)=1$ and
  $d(\underline{m}-e_2)=0$, then one has only $d(\underline{m}-\underline{1})=1$ and
  $c(\underline{m})=1-1-0+1=1$, hence $\underline{m}-\underline{1}$
  cannot be maximal;
  \item[--] if $d(\underline{m}-e_1)=0$ and
  $d(\underline{m}-e_2)=1$, then one has again $d(\underline{m}-\underline{1})=1$ and
  $c(\underline{m})=1$, where again $\underline{m}-\underline{1}$ has to be
  non--maximal;
  \item[--] if $d(\underline{m}-e_1)=0=d(\underline{m}-e_2)=0$, then
  $d(\underline{m}-\underline{1})=0$ and $c(\underline{m})=1$.
\end{itemize}

On the other hand, if $\underline{m}$ is non--maximal, it holds
that
\begin{itemize}
  \item[--] if $d(\underline{m}-e_1)=0$ and $d(\underline{m}-e_2)=1$,
  then $d(\underline{m}-\underline{1})=0$ and the dimension $c(\underline{m})$ 
  has to be equal to $c(\underline{m})=1-0-1+0=0$.
  \item[--] if $d(\underline{m}-e_1)=0$ and
  $d(\underline{m}-e_2)=2$, then we have only $d(\underline{m}-\underline{1})=1$ and
  $c(\underline{m})=1-0-2+1=0$;
  \item[--] if $d(\underline{m}-e_1)=1=d(\underline{m}-e_2)=1$, then again
$d(\underline{m}-\underline{1})=1$ and analogously 
  $c(\underline{m})=1-1-1+1=0$;
  \item[--] if $d(\underline{m}-e_1)=1$ and $d(\underline{m}-e_2)=2$,
  then one has either (i) $d(\underline{m}-\underline{1})=1$, therefore
$c(\underline{m})=1-1-2+1=-1$ and $\underline{m}-\underline{1}$ has to be
maximal by Lemma \ref{lem:propertiessemigroup}(2); or (ii)
$d(\underline{m}-\underline{1})=2$, hence
$c(\underline{m})=1-1-2+2=0$.
\end{itemize}

The remaining case is $d(\underline{m})=0$, which is investigated in a similar manner:

\begin{itemize}
  \item[--] if $d(\underline{m}-e_1)=d(\underline{m}-e_2)=0$, then
  $d(\underline{m}-\underline{1})=0$ and $c(\underline{m})=0$;
  \item[--] if $d(\underline{m}-e_1)=0$ and $d(\underline{m}-e_2)=1$,
  then $d(\underline{m}-\underline{1})=0$ and
  $c(\underline{m})=0$;
  \item[--] if $d(\underline{m}-e_1)=1$ and $d(\underline{m}-e_2)=0$
  one has again $d(\underline{m}-\underline{1})=0$ and
  $c(\underline{m})=0$;
  \item[--] if $d(\underline{m}-e_1)=d(\underline{m}-e_2)=1$,
  then $d(\underline{m}-\underline{1})=2$ and
  $c(\underline{m})=0$. 
\end{itemize}
This finishes the proof. \qed
\medskip

\begin{corollary}
\[
L_{\Gamma_{P_1,P_2}}(\underline{t})=(1-t_1 \cdot t_2)
\sum_{\underline{m} \in \mathcal{M}_{P_1,P_2}}
\underline{t}^{\underline{m}}.
\]
\end{corollary}

\dem~ The basic idea of the proof is to observe
\begin{eqnarray}
L_{\Gamma_{P_1,P_2}}(\underline{t})& = & \sum_{\underline{m} \in
\mathcal{M}_{P_1,P_2}} \underline{t}^{\underline{m}}  +
\sum_{\underline{m}-\underline{1} \in \mathcal{M}_{P_1,P_2}}
\underline{t}^{\underline{m}}\nonumber \\
& = & \sum_{\underline{m} \in \mathcal{M}_{P_1,P_2}}
\underline{t}^{\underline{m}} - \sum_{\underline{m} \in
\mathcal{M}_{P_1,P_2}}
\underline{t}^{\underline{m}+\underline{1}}. \nonumber
\end{eqnarray}
\qed

\begin{theorem}
Let $\mathcal{C}:=\mathcal{C}_{P_1,P_2}$ be the fundamental corner
of the semigroup $\Gamma_{P_1,P_2}$ and consider the set
$\mathcal{M}_{\mathcal{C}}:=\mathcal{M}_{P_1,P_2} \cap
\mathcal{C}$. Then we have
\[
L_{\Gamma_{P_1,P_2}}(\underline{t})=(1-t_1 \cdot t_2)
\sum_{\underline{m} \in \mathcal{M}_{\mathcal{C}}}
\underline{t}^{\underline{m}}.
\]
\end{theorem}

\dem~ We denote
$\underline{\vartheta}^{\prime}:=(\vartheta,-\vartheta)$. The
series $L_{\Gamma_{P_1,P_2}}(\underline{t})$ can be written as
\begin{align*}
L_{\Gamma_{P_1,P_2}}(\underline{t}) & = (1-t_1 \cdot t_2)
\cdot\sum_{\underline{m} \in \mathcal{M}_{P_1,P_2}}
\underline{t}^{\underline{m}} \\
&= (1-t_1 \cdot t_2) \cdot \big ( \sum_{\underline{m}^{\prime} \in
\mathcal{M}_{\mathcal{C}}} \underline{t}^{\underline{m}^{\prime}}
+ \sum_{\substack{\underline{m}^{\prime} + \lambda
\underline{\vartheta}^{\prime}
\\\lambda \in \mathbb{Z}_{\ge 0} \setminus \{0\}}}
\underline{t}^{\underline{m}^{\prime}+\lambda
\underline{\vartheta}^{\prime}} +
\sum_{\substack{\underline{m}^{\prime} - \lambda
\underline{\vartheta}^{\prime}
\\ \lambda \in \mathbb{Z}_{\ge 0} \setminus \{0\}}}
\underline{t}^{\underline{m}^{\prime}-\lambda
\underline{\vartheta}^{\prime}} \big ).
\end{align*}
The proof is completed by showing that the latter two summands are
$0$. \qed

\begin{corollary}
We have
\begin{align*}
P_{\Gamma_{P_1,P_2}}(\underline{t}) & =  \frac{L_{\Gamma_{P_1,P_2}}(\underline{t})}{(1-t_1)(1-t_2)} \\
 & = \frac{1-t_1 t_2}{(1-t_1)(1-t_2)} \sum_{\underline{m} \in \mathcal{M}_{\mathcal{C}}}
 \underline{t}^{\underline{m}}.
\end{align*}
\end{corollary}

\section{Symmetry and functional equations} \label{sec:3}

\subsection{The symmetry of the semigroup and affine embeddings}

Some important considerations about the symmetry of the
Weierstra\ss~semigroup at several points were introduced by
Delgado in \cite{delgado}, being the most remarkable the following
result.

\begin{theorem}[Delgado] \label{thm:defelix}
Let $X$ be a reducible projective space curve of arithmetical
genus $g_X$ and $X^{\prime}=X \setminus \{P_1, \ldots, P_r\}$, where
$P_i$ are smooth points for all $i \in \{1, \ldots, r\}$. Then, the
following statements are equivalent:
\begin{enumerate}
    \item $X^{\prime}$ is an affine complete intersection.
    \item There exists a canonical divisor $K$ such that
    \[
    \mathrm{supp}(K) \subset \{P_1, \ldots, P_r\}.
    \]
    \item There exists a relative maximal $\underline{\tau}=(\tau_1, \ldots,
    \tau_r)$ in the semigroup $\Gamma_{\underline{P}}$ such that $\sum_{i=1}^{r} \tau_i = 2g_X - 2
    + r$.
    \item There exists $\underline{\sigma} = (\sigma_1, \ldots, \sigma_r) \in
    \Gamma_{\underline{P}}$ such that $\sum_{i=1}^{r} \sigma_i= 2g_X-2+r$ and
    $\Gamma_{\underline{P}}$ is symmetric with respect to $\underline{\sigma}$.
\end{enumerate}
\end{theorem}

\dem~ The equivalences $(2) \Leftrightarrow (3) \Leftrightarrow
(4)$ go straight as in \cite{delgado}. The rest of the subsection
is devoted to present the suitable results to reach the
equivalence $(1) \Leftrightarrow (2)$. \qed
\medskip

The rest of the Subsection is devoted to present in a 
self--contained way the proof ``$(1) \Leftrightarrow (2)$'' of 
Theorem \ref{thm:defelix}. (A clear ordering in the steps is missing in the 
proof of Delgado \cite[p.~630]{delgado}). First, we collect now some known results 
due to Serre (cf. \cite{serre3}).

\begin{lemma}  \label{cor:prop2}
Let $R$ be a ring. Assume that every projective $R$--module of
ranks $1$ and $2$ are free. Let $I$ be a nonzero ideal of $R$ with
projective dimension less than or equal to $1$. The following
statements are equivalent:
\begin{enumerate}
    \item $I$ may be generated by two elements.
    \item $\ext_R^1(I,R)$ is a cyclic $R$--module.
\end{enumerate}
\end{lemma}

\dem~ See \cite[\S 2.4.,~Corollaire]{serre3}. \qed
\medskip

Let $V$ be a non--singular variety of dimension $r$,
$\mathcal{O}_V$ be the affine coordinate ring of $V$ and
$\Omega_V$ be the sheaf of differential forms of degree $r$ over
$V$. It is a locally free sheaf of rank $1$. Let $W$ be a
subvariety of $V$ of codimension $h$.
\begin{definition}
The \emph{module of differential forms} on $W$ is defined to be
\[
\Omega_W := \ext_{\mathcal{O}_{V}}^h(\mathcal{O}_{W},\Omega_V).
\]
\end{definition}

\begin{remark} \label{note:vcbundle}
For any ring $R$, there is a correspondence between free finitely
generated projective $R$--modules of finite rank and trivial
vector bundles on $\mathrm{Spec}(R)$ of the same rank.
\end{remark}

The next result was also proven by Serre (cf. \cite[\S
2.7.,~Prop.~6]{serre3}). By the sack of completeness, we rewrite
his proof, adding further details.

\begin{proposition}\label{prop:importante}
Let $V$ be a non--singular affine variety over which every vector
bundle of rank $1$ is trivial. Let $W \subset V$ be a 
Cohen--Macaulay subvariety of $V$ of codimension $2$. Then, if $W$ is a
complete intersection on $V$, then $\Omega_W \cong \mathcal{O}_W$.
Conversely, assuming moreover that every vector bundle of rank $2$
on $V$ is trivial, we have that $W$ is a complete intersection on
$V$ if $\Omega_W \cong \mathcal{O}_W$.
\end{proposition}

\dem~The sheaf $\Omega_V$ is locally free of rank $1$. Every
algebraic vector bundle of rank $1$ is trivial by hypothesis, then
$\Omega_V \cong \mathcal{O}_V$. Furthermore, if $W$ is a complete
intersection, we have the exact sequence
\[
0 \rightarrow \mathcal{A}_W \rightarrow \mathcal{O}_V \rightarrow
\mathcal{O}_W \rightarrow 0,
\]
where $\mathcal{A}_W$ denotes the coherent ideal sheaf defining
$W$. By applying the functor
$\mathrm{Ext}_{\mathcal{O}_V}(-,\mathcal{O}_V)$ to this exact
sequence, the necessity is proven.

Conversely, let $\mathcal{O}_V$ be the affine coordinate ring of
$V$, and $\mathcal{A}_W$ the ideal defining $W$. Then the sheaf
associated to $\mathrm{Ext}_{\mathcal{O}_V}^2 \left
(\mathcal{O}_W, \mathcal{O}_V \right ) \cong
\mathrm{Ext}_{\mathcal{O}_V}^2 \left (\mathcal{O}_W, \Omega_V
\right )$ is isomorphic to $\Omega_W$, and $\Omega_W \cong
\mathcal{O}_W$, thus 
\[
\mathrm{Ext}_{\mathcal{O}_V}^2 \left
(\mathcal{O}_W,\mathcal{O}_V \right ) \cong \mathcal{O}_W
\]
and so
it is cyclic. On the other hand, since $\mathcal{O}_W$ is locally
Cohen--Macaulay, the homological dimension of the ideal
defining $W$ is at most $1$. Remark \ref{note:vcbundle} allows us to 
apply Lemma \ref{cor:prop2} and see that this ideal may be
generated by two elements, i.e., $W$ is a complete intersection on
$V$. \qed

\begin{remark} \label{note:quillen}
Let $R$ be a principal ideal domain. Then every finitely generated
projective $R[x_1, \ldots , x_n]$--module is free (see
\cite[Theorem 4, p.~169]{quillen}).
\end{remark}

Finally, we reach the statement $(1) \Leftrightarrow (2)$ of
Theorem \ref{thm:defelix}.

\begin{corollary} \label{cor:final}
Every non--singular affine curve is a complete intersection if and
only if there exists a canonical divisor with support totally
contained in the places at infinity of the curve.
\end{corollary}

\dem ~ The existence of a canonical divisor totally supported in
places at infinity for the curve is equivalent to the existence of
a trivial canonical line bundle. Let $\mathbb{A}^n =
\mathrm{Spec}(\ff[x_1, \ldots , x_n])$. Since $\ff$ is a field,
every finitely generated projective module over $\ff[x_1, \ldots,
x_n]$ is free by Remark \ref{note:quillen}, and that means
by Remark \ref{note:vcbundle} that every vector bundle on
$\mathrm{Spec}(\ff[x_1,\ldots,x_n])=\mathbb{A}^n$ is trivial. Let
$W$ be a smooth affine curve, then it is Cohen--Macaulay. Now,
from Proposition~\ref{prop:importante}, taking $V=\mathbb{A}^n$
and $W$ the smooth affine curve, that is the case if and only if
$W$ is a complete intersection. \qed

\subsection{Functional equations}

Poincar\'e series associated with value se\-migroups of curve singularities, 
or zeta functions associated with local singular rings,  satisfy functional 
equations whenever the base ring is Gorenstein  (see \cite{de2} and \cite{stohr}; also \cite{mozu}). 
The same holds in our global context for complete intersections.

\begin{definition}
The point $\underline{\sigma}$ of the last assertion of Theorem
\ref{thm:defelix} is called the \emph{symmetry point} of
$\Gamma_{\underline{P}}$.
\end{definition}

For the case $r=1$, the symmetry point of $\Gamma_P$ is
the conductor, which is equal to $c(\Gamma_P)=2g_X$ (cf. \cite{sathaye}). Then
we can deduce easily the functional equations.

\begin{lemma}
Let $X$ be a curve. Let $P$ be a rational point of $X$ such that
$X \setminus P$ is an affine complete intersection. We have
\begin{align*}
L_{\Gamma_P}(t) &=  - t^{2g_X} L_{\Gamma_P}(t^{-1}),\\
P_{\Gamma_P}(t)  &= t^{2g_X-1} P_{\Gamma_P} (t^{-1}).
\end{align*}
\end{lemma}

\dem~ Let $n \in \mathbb{Z}$. By the definition of conductor of a
numerical semigroup, we know that $n \in \Gamma_P$ if and only if
$2g_X-1-n \notin \Gamma_P$; it implies that $d(2g_X-1-n)=0$ if and
only if $d(n)=1$. It is easily checked that
\[
t^{2g_X-1} P_{\Gamma_P}(t^{-1})  = \sum_{m \in \mathbb{Z}} d (m)
t^{2g_X-1-n}
 =  \sum_{n \in \Gamma_P} t^n
 =  P_{\Gamma_P} (t).
\]
The functional equation for $L_{\Gamma_P}(t)$ follows from a
simple computation. \qed
\medskip

The next result establishes the functional equation for the Poincar\'e
series associated with two points.

\begin{proposition}
Let $X$ be a curve, let $P_1,P_2$ be two rational points on $X$ such that $X \setminus \{P_1,P_2\}$ is an affine complete intersection. Let us
denote by $\underline{\sigma}$ the symmetry point of the
Weierstra\ss~semigroup $\Gamma_{P_1,P_2}$. We have
\begin{align*}
L_{\Gamma_{P_1,P_2}}(\underline{t})&=
-\underline{t}^{\underline{\sigma}
+ \underline{1}} \cdot L_{\Gamma_{P_1,P_2}}(\underline{t}^{-1} ),\\
P_{\Gamma_{P_1,P_2}}(\underline{t})
&=-\underline{t}^{\underline{\sigma} } \cdot
P_{\Gamma_{P_1,P_2}}(\underline{t}^{-1} ).
\end{align*}
\end{proposition}

\dem~ By Proposition \ref{prop:importante}, we have
\[
L_{\Gamma_{P_1,P_2}}(\underline{t})= (1-t_1 \cdot t_2) \cdot
\sum_{\underline{m} \in \mathcal{M}_{P_1,P_2} }
\underline{t}^{\underline{m}}.
\]

Moreover, for a maximal point $\underline{m} \in \mathcal{M}_{P_1,P_2}$, there
exists $\underline{n} \in \mathcal{M}_{P_1,P_2}$ such that
$\underline{n}+\underline{m}=\underline{\sigma}$. Hence
\begin{align*}
  L_{\Gamma_{P_1,P_2}}(\underline{t})& =  (1-t_1 \cdot t_2)\cdot \sum_{\underline{m} \in \mathcal{M}_{P_1,P_2} }
\underline{t}^{\underline{m}} \\
   & =  (1-t_1 \cdot t_2) \cdot \sum_{\underline{n} \in \mathcal{M}_{P_1,P_2} }
\underline{t}^{\underline{\sigma}-\underline{n}} \\
   & = -\underline{t}^{\underline{\sigma} + \underline{1}} \cdot
   L(\underline{t}^{-1}).
\end{align*}
On the other hand, by the definition and the functional equation
of $L_{\Gamma_{P_1,P_2}}(\underline{t})$, we have
\begin{align*}
  P_{\Gamma_{P_1,P_2}} \big (\underline{t}^{-1} \big )& =  \frac{L_{\Gamma_{P_1,P_2}} \big (\underline{t}^{-1} \big )}{(1-t_1^{-1})(1-t_2^{-1})}\\
   & =  \frac{-\underline{t}^{-\underline{\sigma}} \cdot
   L_{\Gamma_{P_1,P_2}}(\underline{t})}{(1-t_1)(1-t_2)}
   \\
   & =  - \underline{t}^{-\underline{\sigma}} \cdot
   P_{\Gamma_{P_1,P_2}}(\underline{t}).
\end{align*} \qed

\section*{Acknowledgements}

The author is thankful to Professor F\'elix Delgado de la Mata for stimulating conversations on the combinatorial structure of value semigroups, as well as to the
referee of this paper for his/her helpful comments.

\end{document}